\documentclass[12pt]{article}
\usepackage{graphicx}
\usepackage{amsmath}
\usepackage{amsfonts}
\usepackage{amsmath,amsfonts,amssymb}
\usepackage[pdftex,bookmarks,colorlinks]{hyperref}
\newcounter{alphthm}
\newenvironment{proof}{\par\noindent\textbf{Proof}.}{\hfill$\square$}
\newtheorem{prop}{Proposition}[section]
\newtheorem{thm}{Theorem}[section]
\newtheorem{lem}[thm]{Lemma}

\newtheorem{cor}{Corollary}[section]

\newtheorem{definition}{Definition}[section]
\newtheorem{remark}{Remark}[section]
\newtheorem{example}{Example}

\newcommand{\bpf}{\begin{proof}}
\newcommand{\epf}{\end{proof}}
\newcommand{\bl}{\begin{lem}}
\newcommand{\el}{\end{lem}}
\newcommand{\bp}{\begin{prop}}
\newcommand{\ep}{\end{prop}}
\newcommand{\bt}{\begin{thm}}
\newcommand{\et}{\end{thm}}
\newcommand{\bc}{\begin{cor}}
\newcommand{\ec}{\end{cor}}
\newcommand{\R}{I\!\! R}
\def\nn{\nonumber}

\title{On Sobolev spaces and density theorems on Finsler manifolds}
\author{B. Bidabad\thanks{The corresponding author, bidabad@aut.ac.ir} and A. Shahi}
\date{}
\def\beq{\begin{equation}}
\def\eeq{\end{equation}}
\def\nn{\nonumber}
\begin{document}
\maketitle
\begin{abstract}
Here, a natural extension of Sobolev spaces is defined for a  Finsler structure $F$ and it  is shown that the set of all real $C^{\infty}$ functions  with compact support on a forward geodesically complete Finsler manifold $(M, F)$, is dense in the extended Sobolev space $H_1^p (M)$.
As a consequence,  the weak solutions $u$ of the Dirichlet equation $\Delta u=f$ can be approximated by $C^\infty$ functions with compact support on $M$.
Moreover, let $W \subset M$ be a regular domain with the $C^r$ boundary $\partial W$, then the set of all real functions in $C^r (W) \cap C^0 (\overline W)$ is dense in $H_k^p (W)$, where $k\leq r$. Finally, several examples are illustrated and   sharpness of the inequality $k\leq r$ is shown.
\end{abstract}
\vspace{1cm}

{\footnotesize\textbf{Keywords}: {  Density theorem, Sobolev spaces, Dirichlet problem, Finsler space.}}

{\footnotesize\textbf{AMS subject classification}: {53C60, 46E35}}
\section{Introduction}
A Sobolev space is a vector space of functions endowed with a norm which is a combination of $L^p -$norm of the function itself and its derivatives up to a certain order.
Its objectives are to deal with some basic PDE problems on manifolds equipped with Riemannian metrics. For instance, the Yamabe problem asserts that for any compact Riemannian manifold $(M,g)$ of dimension $n\geq 3$, there always exists a metric with constant scalar curvature. Clearly, solutions of Yamabe elliptic equation on Riemannian manifolds are laying in the Sobolev space $H_1 ^2 (M)$, see for instance \cite{Au2}. Currently, the question of elliptic equations on some natural extensions of Riemannian spaces, particularly on Finsler manifolds are extensively studied, see for instance \cite{AR, BSY, JS, KR, La}.
Another natural question is to wonder whether a function can be approximated by another one with better properties.
Density problems permit to investigate this question and find conditions under which a function on a Sobolev space can be approximated by smooth functions with compact support.

 Historically, one of the significant density theorems is proved by S. B. Myers \cite{My} in 1954 for compact Riemannian manifolds and then in 1959 by M. Nakai for finite-dimensional Riemannian manifolds.
  Next, in 1976, T. Aubin has investigated density theorems on Riemannian manifolds, cf. \cite{Au2}.
Y. Ge and Z. Shen \cite{GS} defined a canonical energy functional on Sobolev spaces and investigated the eigenvalues and eigenfunctions related to this functional, on a compact reversible Finsler manifold. Next, Y. Yang has defined a Sobolev space on a reversible Finsler manifold $(M, F)$ by using the osculating Riemannian metric and the corresponding Levi-Civita connection on the underlying manifold $M$,  cf., \cite{Ya}. In 2011 the Myers-Nakai theorem is extended to the Finsler manifolds of class $C^k$, where $k \in N\cup \{\infty \} $, cf., \cite{JS}. Lately, S. Ohta has studied many Sobolev inequalities in Finsler spaces, see \cite{Oh}.

Recently the present authors have studied some natural extensions of Riemannian results, more or less linked to this question see for instance \cite{Bi, BS2, BY}.

  In the present work, a natural extension of Riemannian metrics is considered, and solutions to the above questions are studied. More intuitively,  a Sobolev space is defined by considering a Riemannian metric on the sphere bundle $SM$, rather than the manifold $M$, and following Theorems are proved. 
    Denote by $D(M)$ the set of all real $C^\infty$ functions with compact support on $M$ and let $ \mathop {H^p_1}\limits^{\circ~} (M)$ be the closure of $D(M)$ in $H^p_1 (M),$ where  $(M,F)$ is an $n-$dimensional $C^\infty$ Finsler manifold, $p\geq 1$ a real number, $k$ a non-negative integer and $H_k^p (M)$ certain Sobolev space determined by the Finsler structure $F$.
 \bt \label{Th1}
Let $(M,F)$ be a forward (or backward) geodesically complete Finsler manifold, then $\mathop{H^p_1}\limits^{\circ~} (M)=H_1^p (M)$.
\et
As an application for a real function $f:M\longrightarrow{\R}$ on a compact $C^\infty,$ reversible Finsler manifold for which $\int _M fdv_F =0$, the weak solutions $u$ of the Dirichlet equation $\Delta u=f$ can be approximated by $C^\infty$ functions with compact support on $M$.

We provide also some examples.

 Let $W \subset M$ be a s-dimensional regular domain with $C^r$ boundary $\partial W$, then $(\overline W,F)$ is called a Finsler manifold with $C^r$ boundary.\\

 \bt \label{Th2}
Let $(\overline W,F)$ be a compact Finsler manifold with $C^r$ boundary. Then $C^r (\overline W)$ is dense in $H_k^p (W)$, for $k \leq r.$
\et
Next, using a counter example we show that the assumption $k \leq r$, in  Theorem \ref{Th2} is sharp and can not be omitted.
As a consequence of the above density theorems, we can approximate solutions of partial differential equations on a Sobolev space determined by $F$, with $C^\infty$ or $C^r$ functions on $(M, F)$ and hence study some recent problems on Finsler geometry, for instance, Ricci flow, Yamabe flow, etc.

It should be recalled that the new definition of Sobolev space in Finsler geometry introduced in the present work, is a more general definition and has the following advantages.
\begin{itemize}
\item
 This definition of Sobolev space, reduces to that of Ge and Shen, in the case $k=1$ and $p=2$, provided the underlying manifold is closed and the Finsler structure is reversible, see Remark \ref{Rem; Shen's definition}.
\item In this approach, the reversibility condition on the Finslerian structure is not required.
\item The present definition is applicable also to the geometric objects defined on $TM$.
\item This approach makes possible to study the Sobolev norms of horizontal curvature tensor and its covariant derivatives up to the certain order $k$.
\item This approach, permits to generalize Theorem\ref{Th1} for $H_k^p (M)$, where $k\geq 2$.
\end{itemize}
 We adopt here the notations and terminologies of \cite{Ak, BCS}, and \cite{Sh} and recall that all the Finsler manifolds in the present work are assumed to satisfy in Remark \ref{stry}.
\section{Preliminaries and terminologies}

Let $M$ be a connected differentiable manifold of dimension $n$. Denote the bundle of tangent vectors of $M$ by
$\pi_1:TM\longrightarrow M$, the fiber bundle of non-zero tangent
vectors of $M$ by $\pi:TM_0\longrightarrow M$ and the
pulled-back tangent bundle by $\pi^*TM$. A point of $TM_0$ is denoted
by $z=(x,y)$, where $x=\pi z\in M$ and $y\in T_{\pi z}M$. Let $(x,U)$ be a local chart on $ M$ and $(x^i,y^i)$ the induced local coordinates on $\pi^{-1}(U)$, where ${\bf y}=y^i\frac{\partial}{\partial x^i}\in T_{\pi z}M$, and $i$ running over the range $1,2,...,n$. A (globally defined) \emph{Finsler structure} on $M$ is a function $F:TM\longrightarrow [0,\infty)$ with the following properties;
$F$ is $C^\infty$ on the entire slit tangent bundle $TM \backslash 0 $;
$F(x,\lambda y )=\lambda F(x,y) ~ \forall \lambda >0$;
the $n\times n$ Hessian matrix $(g_{ij}) =\frac {1}{2} ([F^2]_{y^i y^j})$ is positive-definite at every point of $TM_0$.
The pair $(M,F)$ is called a {\it Finsler manifold}.
Given the induced coordinates $(x^i,y^i)$ on $TM$, coefficients of spray vector field are defined by $G^{i}=1/4 g^{ih}(\frac{\partial^{2}F^{2}}{\partial y^{h}\partial x^{j}}y^{j}-\frac{\partial F^{2}}{\partial x^{h}})$.
One can observe that the pair $\{\delta/\delta x^i,\partial/\partial y^i \}$ forms a horizontal and vertical frame for $TTM$, where ${\frac {\delta}{\delta x^i}}:={\frac {\partial}{\partial x^i}}-G_i^j{\frac {\partial}{\partial y^j}}$, $G_i^j :={\frac {\partial G^j}{\partial y^i}}$. Denote by $SM$ the sphere bundle, where $SM:=\bigcup\limits _{x\in M} S_xM$ and $S_xM:=\{y\in T_xM | F(y)=1\}$.
The\emph{ Sasakian metric} on $SM$ is defined by
\begin{equation}\label{Sa}
\hat{g}=\delta_{ab}w^a\otimes w^b +\delta _{\alpha \beta} w^{n+\alpha}\otimes w^{n+\beta},
\end{equation}
where $a,b=1,..,n$ and $\alpha ,\beta =1,...,n-1$, and $\{w^a ,w^{n+\alpha}\}$ is an ordered orthonormal coframe on $SM$, cf., \cite{BL}
The volume element $dV_{SM}$ of $SM$ with respect to the Sasakian metric $\hat{g}$ is
\begin{equation*}
dV_{SM}=w^1 \wedge ... \wedge w^{2n-1}=\sqrt {\det(g_{ij})} dx^1\wedge ...\wedge dx^n \wedge w^{n+1} \wedge ...\wedge w^{2n-1}.
\end{equation*}
This volume element can be rewritten as follows,
\begin{equation}\label{SV}
dV_{SM} =\Omega d\tau \wedge dx,
\end{equation}
where, $ \Omega=\det(\frac{g_{ij}}{F}),~d\tau =\sum\limits_{k=1}^n (-1)^{k-1} y^k dy^1 \wedge ... \hat{dy^k}\wedge ...\wedge dy^n$ and $dx$ is the n-form $ ~dx=dx^1 \wedge ...\wedge dx^n$, cf. \cite{HS}.
We have a volume form on $(M,F)$
\begin{equation}\label{MV}
dV_{F} = (\frac {1}{c_{n-1}} \int_{S_x M} \Omega d\tau ) dx,
\end{equation}
where $c_{n-1}$ denotes the volume of the unit Euclidean sphere $S^{n-1}$, cf. \cite{HS}.

 Let $ \sigma:[a,b]\longrightarrow M$ a piecewise $C^\infty$ curve with the velocity $\frac {d \sigma}{dt} =\frac {d \sigma ^i}{dt} \frac {\partial}{ \partial x^i} \in T_ {\sigma (t)} (M).$
Its integral length is defined by $L(\sigma)= \int _a^b F(\sigma , \frac {d \sigma}{dt}) dt$. For $x_1, x_2 \in M$ denote by $ \Gamma (x_1, x_2)$ the collection of all piecewise $C^\infty$ curves $\sigma:[a,b]\longrightarrow M$ with $\sigma (a)=x_1$ and $\sigma (b)=x_2$ and by $d(x_1, x_2)$ the metric distance from $x_1$ to $x_2$,
\begin{equation} \label{Eq;distance}
d(x_1, x_2)=\inf\limits_{\Gamma (x_1, x_2)} L(\sigma).
\end{equation}
\bl \label{Lemma}
\cite{BCS} Let (M,F) be a Finsler manifold. At any point $x \in M$, there exists a local coordinate system $(\phi ,U)$ such that the closure of U is compact, $\phi :\overline U\longrightarrow {\R}^n$, $\phi (x)=0$ and $\phi$ maps $U$ diffeomorphically onto an open ball of ${I\!\!R}^n.$
\el
A Finsler manifold is said to be forward (resp. backward) geodesically complete if every geodesic $\gamma (t),~ a \leq t <b,$ parameterized to have constant Finslerian speed, can be extended to a geodesic on $a \leq t<\infty.$ (resp. $-\infty <t\leq b$). If the Finsler structure $F$ is reversible, then $d$ is symmetric. In this case, forward completeness is equivalent to backward completeness.
Compact Finsler manifolds at the same time are both forward and backward complete, whether $d$ is symmetric or not.
\section{A Sobolev space on Finsler manifolds}

Let $(M,F)$ be a $C^\infty$ Finsler manifold. For any real function $u$ on $M$ we denote again $u\circ\pi$ by $u$.
 The $j$th covariant derivative of $u$ is denoted by $\nabla ^j u$, where $\nabla$ is horizontal covariant derivative of Cartan connection, $j$ is a nonnegative integer, hence $\nabla ^0 u=u$. Let us denote the inner scalar product on $SM$ with respect to the Sasakian metric (\ref{Sa}) by $(.,.)$ and $|\nabla ^j u|^2 =(\nabla ^j u ,\nabla ^j u)$. We denote by $C^p _k (M)$ the space of smooth functions $u \in C^ \infty (M)$ such that $|\nabla ^j u| \in L^p (SM)$ for any $j$ run over the range $0,1, ...,k$ and $p\geq 1$, that is
$$C^p _k (M)=\{u\in C^ \infty (M) : \forall j=0,1,...,k,~ \int_{SM} [(\nabla ^j u ,\nabla ^j u )] ^{\frac{p}{2}}dV_{SM}< \infty \}.$$

\begin{remark}\label{stry}
It is well known that $S^{n-1}$ is diffeomorphic  to the $S_x M$, where $x\in M$. Let $A$ be the radial projection from $S_x M =\{(y^i) \in {\R}^n: F(x, y^i \frac{\partial}{\partial x^i})=1 \} \subset {\R}^n$ onto the unit sphere $S^{n-1} \subset {\R}^n$ and $(\det J(A))$ determinant of its Jacobian.
 Everywhere in this paper we assume that there exists a positive real  number $R>0$ such that $(\det J(A)) \sqrt{\det g_{ij}}\geq \frac{R}{c_{n-1}}$, where ${c_{n-1}}$ is the volume of $S^{n-1}$. Due to the dependance of Cartan covariant derivatives to the direction on Finsler cases we have to consider the above inequality. Note that every compact Finsler manifold satisfies the preceding inequality, due to the compactness of $SM$.
\end{remark}
\bp \cite{VW} \label {TTT}
Let $(x,\Omega)$ be a local coordinate chart on $M$ and $f:SM \subset TM_0 \rightarrow {\R}$ an integrable function with support in $\pi ^{-1} (\Omega)$. Then we have
\begin{equation} \label{prop}
\int _{SM} f(x,y)dV_{SM} =\int _{\Omega} (\int _{S^{n-1}} f(x,\frac{y}{F})\frac{det(g_{ij})}{F^n}d\sigma)dx,
\end{equation}
where $d\sigma$ is the standard volume form on $S^{n-1}$, $dx=dx^1 \wedge ....\wedge dx^n$, $y=y(x,\theta)$
for $(x,\theta) \in \Omega \times S^{n-1}$ and $\theta =(\theta ^1 ,...,\theta ^n)$ are local coordinate on $S^{n-1}$.
\ep
\begin{definition}\label{Df}
The Sobolev space $H_k ^p (M)$ is the completion of $C_k ^p (M)$ with respect to the norm
$$\parallel u\parallel _{H_k ^p (M)}=\sum\limits _{j=0} ^k (\int _{SM} [(\nabla ^j u ,\nabla ^j u)]^{\frac{p}{2}}dV_{SM})^{\frac {1}{p}},$$
where $p \geq 1$ is a real number.
\end{definition}
Let $f:M\rightarrow {\R}$ be a real function, then using the volume form $dV_F$ defined by (\ref{MV}) we can consider the space of $L^p-$norm as follows
\begin{equation}
L^p (M)=\{f:M\rightarrow {\R}~~ \textit{is measurable} : \int _M f^p dV_F <\infty \}.
\end{equation}

 \bl \label{SSSS}
The Finslerian Sobolev space $H_k ^p (M)$ is a subspace of $L^p (M)$.
\el
\begin{proof}
For all $u \in H_k ^p (M)$,  by Definition \ref{Df} we have $\int_{SM} |u|^p dV_{SM} <\infty$.
An appropriate choice of $\{\Omega _i\}_{i=1}^\infty$, together with Proposition \ref{TTT} and relations (\ref{SV}) and (\ref{MV}) leads to
\begin{align*}
\int _{SM} |u|^p dV_{SM} &= \sum\limits_{i=1}^\infty \int _{\pi ^{-1}(\Omega _i)} |u|^p dV_{SM}\\
&=\sum\limits_{i=1}^\infty \int _{\Omega _i} \int _{S^{n-1}} |u|^p \frac{det(g_{ij})}{F^n} d\sigma dx\\
&=\sum\limits_{i=1}^\infty \int _{\Omega _i} \int _{S_x M} |u|^p \frac{det(g_{ij})}{F^n} (\det J(A))~ dV_{S_x M} dx\\
&=\sum\limits_{i=1}^\infty \int _{\Omega _i} \int _{S_x M} |u|^p \frac{det(g_{ij})}{F^n} (\det J(A)) ~\sqrt{\det g_{ij}} d\tau dx.
\end{align*}
By assumption there exists $R>0$ such that $(\det J(A)) ~\sqrt{\det g_{ij}} \geq\frac{R}{c_{n-1}}$, hence
\begin{align*}
\sum\limits_{i=1}^\infty \int _{\Omega _i} \int _{S_x M} |u|^p \frac{det(g_{ij})}{F^n} &(\det J(A)) ~\sqrt{\det g_{ij}} d\tau dx\\
&\geq R \sum\limits_{i=1}^\infty \int _{\Omega _i} |u|^p (\frac{1}{c_{n-1}}\int _{S_x M}  \Omega d\tau ) dx\\
&= R \sum\limits_{i=1}^\infty \int _{\Omega _i} |u|^p dV_F\\
&=R \int _M |u|^p dV_F.
\end{align*}
 Therefore $u \in L^p (M)$ and proof is complete.
\end{proof}

The following example shows that  the assumption $(\det J(A)) ~\sqrt{\det g_{ij}} \geq\frac{R}{c_{n-1}}$ is essential and can not be dropped.
\begin{example}
Let $M={\R}^2$, hence $SM\simeq {\R}^2 \times S^1$. Consider the function $f:{\R}^2 \rightarrow {\R}$ defined by
\begin{equation*}
f(x,y)=\left\{ \begin{array}{ll}1& 0<y<1\\ 0& o. w.\end{array}.\right.
\end{equation*}
Choose a metric on $M$, such that the fibers  $S_{(x,y)} M$ of $SM$ have the radius $r_x=e^{-x^2}$. Let $U=\R\times (0,1)$ which leads to
\begin{align*}
\int _{SM} f dV_{SM} &=\int_{U\times S_{(x,y)}M} d\theta dx dy= \int_U (\int _{S_{(x,y)}M} d\theta) dxdy\\
&=\int _0 ^1 \int _{-\infty} ^{+\infty} (2 \pi e^{-x^2})dx dy <\infty .
\end{align*}
On the other hand
\begin{align*}
\int _{M} f d\mu =\int _0 ^1 \int _{-\infty} ^{+\infty} 1dx dy =\infty .
\end{align*}
Therefore $f\notin L^1 ({\R}^2)$.
\end{example}
\begin{remark}\label{Rem;Shen's definition}
Ge and Shen  in \cite{GS}, defined a  norm of  Sobolev spaces on a closed reversible  Finsler manifold, for $k=1$ and $p=2$ as follows
\begin{align} \label{GS Def}
\parallel u \parallel _{_{GS}}= (\int _M u^2 dV_F)^{\frac{1}{2}} + (\int _M (F(\nabla u)) ^2 dV_F)^{\frac{1}{2}},
\end{align}
where $(\nabla u) (x)$ is the gradient of $u$ in $x\in M$. A similar argument as  in the proof of Lemma \ref{SSSS}, shows that the definition of  Sobolev space given in the present work, reduces to that of \cite{GS}, for $k=1$ and $p=2$, on a closed reversible  Finsler manifold.
\end{remark}
\begin{remark}\label{Rem;3}
In Definition \ref{Df} we use an inner product on $SM$ to define a Sobolev space and naturally it has a structure of  vector space, while the definition given in \cite{GS} is not a vector space on complete Finsler manifolds in general. In fact,  Krist\'{a}ly and  Rudas show that on $(B^n (1),F)$, where $B^n (1)$ is an $n-$dimensional unit ball of ${\R}^n$ and $F$ is a Funk metric, the Sobolev space defined in \cite{GS} has no more structure of a vector space, due to the irreversibility of $F$, cf., \cite{KR}.
\end{remark}

 Let $J$ be a nonnegative, real-valued function, in the space of $C^\infty $ functions with compact support on ${\R}^n,$ denoted by $C_0^\infty ({\R}^n)$ and having properties :\\
(i)$ J(x)=0$ if $| x| \geq 1 .$\\
(ii)$ \int _{{\R}^n} J(x)dx =1.$\\
Consider the function $J_\epsilon (x)=\epsilon ^{-n} J(\frac{x}{\epsilon})$ which is nonnegative in $C_0^\infty ({\R}^n)$ and satisfies :\\
(i) $J_\epsilon (x) =0$ if $| x| \geq \epsilon ,$\\
(ii) $ \int _{{\R}^n} J_\epsilon (x)dx =1.$\\
$J_\epsilon$ is called a \emph{mollifier} and the convolution $J_\epsilon * u(x) := \int _{{\R}^n} J_\epsilon (x-y) u(y)dy$, defined for the function $u$ is called a \emph{regularization} of $u.$

 \bl\cite{Ad} \label{lem1}
Let u be a function defined on $\Omega \subset {\R}^n$ and vanishes identically outside the domain $\Omega :$\\
(a) If $u \in L^1_{Loc}( \overline {\Omega }) $ then $J_\epsilon * u \in C_0^\infty ({\R}^n).$\\
(b) If  $supp(u)\subset \Omega,$ then $J_\epsilon * u \in C^\infty _0 (\Omega)$ provided $\epsilon < dist(supp(u),\partial \Omega).$\\
(c) If $u \in L^p (\Omega)$ where $1\leq p< \infty$ then $J_\epsilon * u \in L^p (\Omega ).$ Moreover,\\ $\parallel J_\epsilon *u \parallel _p \leq \parallel u \parallel _p$
and $\lim\limits _{\epsilon \rightarrow 0^+} \parallel J_\epsilon *u -u \parallel _p =0$
\el


 \bl \cite{Ad} \label{lem2}
Let $u \in H_k^p (\Omega)$ and $1 \leq p<\infty.$ If $\Omega' \subset \subset \Omega \subset {\R}^n,$ that is $\overline{\Omega'}\subset \Omega$ and $\overline{\Omega'}$ is a compact subset of ${\R}^n$, then $\lim\limits_{\epsilon\rightarrow o^+} J_\epsilon * u=u$ in $H_k^p (\Omega')$.
\el


 \section{Density theorems on Finsler manifolds}
  In the previous section we set necessary tools on $SM$ which permits to use  Aubin's techniques in Finsler geometry, cf. \cite{Au2}.
Let $(M,F)$ be a Finsler manifold and $D(M)$ the space of $C^\infty$ functions with compact support on $M$.
In this section, we use Hopf-Rinow's theorem to introduce the first density theorem on boundaryless Finsler manifolds and investigate another density theorem for Finslerian manifolds with $C^r$ boundary.

 \subsection{Case of manifolds without boundary}~~~~~~~~~~
 \bp \label{fl}
Let $(M,F)$ be a forward geodesically complete Finsler manifold, then any function $\phi \in H_1 ^p (M)$ can be approximated by functions with compact support on $M$.
\ep
\bpf ~Let $\phi \in H_1 ^p (M)$, then by Lemma \ref{SSSS} we have $\phi \in L^p (M)$  and hence it can be approximated by $C^\infty$ functions. Therefore $C^\infty (M)\cap H_1^p (M)$ is dense in $H^p_1 (M)$. To prove Proposition \ref{fl} it suffices to show the assertion for $\phi \in C^\infty (M)\cap H_1^p (M)$.
Let $\phi\in C^\infty (M) \cap H_1^p (M)$ and fix a point $x_0\in M$. By Hopf-Rinow's~theorem for forward geodesically complete Finsler manifolds, every pair of points in $M$, containing $x_0$ can be joined by a minimizing geodesic emanating from $x_0$. Let us consider the Finslerian distance function (\ref{Eq;distance}), and the sequence of functions $\phi_j (x) =\phi (x) f(d(x_0 ,x)-j),$ where $f:{\R}\longrightarrow {\R} $ is defined by
\begin{equation}
f(t)=\left\{ \begin{array}{ll}1& t\leq 0\\ 1-t&0<t<1\\ 0&t\geq 1\end{array}.\right.
\end{equation}
Clearly, $f$ is a continuous and decreasing function and it is differentiable almost everywhere  on ${\R}$. We should prove, each $\phi _j (x)$ has a compact support.
If $d(x_0 ,x) \leq j$ then $d(x_o ,x)-j\leq 0$ or $f(d(x_o ,x)-j)=1$ hence $\phi _j (x)=\phi (x).$
If $d(x_0 ,x)\geq j+1$ then $d(x_o ,x)-j\geq 1$ or $\phi _j (x)=0.$
Thus each $\phi _j$ has a compact support and obviously $\lim\limits _{\scriptscriptstyle j\rightarrow \infty} {\phi _j (x)}=\phi (x)$. This  completes the proof .
\epf
\bp \label{sl}
Let $(M,F)$ be a forward geodesically complete Finsler manifold, then any function $\phi \in H_1 ^p (M)$ can be approximated by functions with compact support in $H_1 ^p (M)$.
\ep
\bpf ~Let $\phi \in H_1 ^p (M) \cap C^\infty (M)$ be an arbitrary function and $\{\phi _j\}_{j=1}^\infty$ the sequence of related functions given in the proof of Proposition \ref{fl}. We show that $\parallel \phi _j -\phi \parallel _{H_1 ^p (M)}$ tends to zero as $j\rightarrow \infty$. Since $\phi$ is a $C^ \infty $ real function on $M$ and $\phi (x)\in H^p_1 (M)$, $(\nabla \phi)_{(x,y)}$ exists and bounded almost everywhere on $(x,y) \in SM$.
It is well known on forward geodesically complete Finsler manifolds, the Finslerian distance function $d$ is $C^\infty$ out of a small neighborhood of $x_0$ and it is $C^1$ in a punctured neighborhood of $x_0$ cf. \cite{Wi}.
By definition of $\phi _j$, $|\phi _j (x)| \leq |\phi (x)|$, where $\phi (x)$ is a $C^\infty$ function lying in $L^p (SM)$. Therefore by the Lebesgue dominated convergence theorem $\phi _j$ and $\phi _j -\phi $ are also in $L^p (SM) $. Recall that
\beq \label{Eq 3.1}
\parallel \phi _j -\phi \parallel _{H_1 ^p (M)}=\sum\limits _{k=0} ^1 (\int _{SM} [(\nabla ^k (\phi _j -\phi) ,\nabla ^k (\phi _j -\phi))]^{\frac{p}{2}}dV_{SM})^{\frac {1}{p}}.
\eeq
We show that equation (\ref{Eq 3.1}) tends to zero as $j$ tends to infinity. First assume support of  $\phi _j -\phi$ is contained in a local chart $(x,\Omega)$ , then by Proposition \ref{TTT} we have
\begin{equation*}
\int _{SM} |\phi _j -\phi|^p dV_{SM} =\int _{\Omega} (\int _{S^{n-1}} |\phi _j -\phi|^p _{(x,\frac{y}{F})}\frac{det(g_{ij})}{F^n}d\sigma)dx.
\end{equation*}
Consider forward metric ball
$B_{x_0}^+ (j)=\{x \in \Omega: d(x_0 ,x)<j\}$, we have
\begin{align} \label{Eq 3.2}
\int_ {SM} |\phi _j -\phi|^p dV_{SM}& \leq \int _{{B_{x_0}^+ (j)}} (\int _{S^{n-1}} |\phi _j -\phi|^p _{(x,\frac{y}{F})}\frac{det(g_{ij})}{F^n}d\sigma)dx \nn\\
&+\int _{\Omega \setminus{B_{x_0}^+ (j)}} (\int _{S^{n-1}} |\phi _j -\phi|^p _{(x,\frac{y}{F})}\frac{det(g_{ij})}{F^n}d\sigma)dx. \nn
\end{align}
When $j$ tends to infinity $\Omega \setminus {B_{x_0}^+ (j)}=\emptyset$, therefore

 \begin{eqnarray*}
\int_ {SM}|\phi _j -\phi|^p dV_{SM}&\leq& \int _{{B_{x_0}^+ (j)}} (\int _{S^{n-1}} |\phi _j -\phi|^p _{(x,\frac{y}{F})}\frac{det(g_{ij})}{F^n}d\sigma)dx.\nn\\
\end{eqnarray*}
On $B_{x_0}^+ (j)$ we have $d(x_0 ,x)<j$ or $\phi_j (x)=\phi (x)$, hence
\begin{equation} \label{r1}
\parallel \phi _j - \phi \parallel _p \leq (\int _{SM}|\phi_j -\phi|^p dV_{SM})^{\frac{1}{p}} \longrightarrow 0.
\end{equation}
Next if support of $\phi _j -\phi$ is not contained in a local chart, then similar to the proof of Lemma \ref{SSSS} by appropriate choice of sequence $\{\Omega _i \}_{i=1} ^\infty$ we can show that $\parallel \phi _j -\phi \parallel _p \rightarrow 0$.
 We prove $\parallel \phi _j - \phi \parallel _{H^p_1 (M)} $ converges to zero. To this end it suffices to show that $\parallel \nabla \phi _j -\nabla \phi \parallel _p $ or $ \parallel \nabla (\phi_j - \phi) \parallel _p $ converges to zero.
By means of Leibnitz's formula for $\phi _j (x)=\phi (x) f(d(x_0, x)-j)$, and triangle inequality we obtain
\begin{equation} \label{lib}
| \nabla \phi _j | \leq | \nabla \phi | + | \phi | \sup\limits _{t \in [0,1]} | f'(t) |.
\end{equation}
Again with Lebesgue dominated theorem, we have $|\nabla \phi _j| \in L^p (SM)$ and hence $\phi _j (x) \in H_1^p (M).$
Repeating above steps for $|\nabla (\phi _j -\phi)|$ instead of $|(\phi _j -\phi)|$ leads to
\begin{equation} \label{r2}
\parallel \nabla \phi _j-\nabla \phi \parallel _p = \parallel \nabla (\phi_j - \phi) \parallel _p \leq(\int _{SM}|\nabla (\phi _j -\phi)|^p dV_{SM} )^{\frac {1}{p}} \longrightarrow 0.
\end{equation}
Therefore by the relations (\ref{r1}) and (\ref{r2}) we obtain
$$\parallel \phi _j -\phi \parallel _ {H_1^p (M)}= \parallel \phi _j - \phi \parallel _p + \parallel \nabla(\phi _j -\phi )\parallel _p \longrightarrow 0.$$ Thus $\phi_j$ converges to $\phi$ in $H_1^p (M)$.

 \epf

 Proof of Theorem \ref{Th1} is an application of Propositions \ref{fl} and  \ref{sl}, similar to that in Riemannian geometry.

 \noindent {\bf Proof of Theorem \ref{Th1}.}
 To prove Theorem \ref{Th1}, by means of Propositions \ref{fl} and \ref{sl}, it remains to approximate each $\phi _j$ by functions in $D(M)$. Let $j$ be a fixed index for which $\phi _j$ has a compact support. Let $K$ be the compact support of $\phi _j$ and $\{V_i \}_{i=1}^m $ a finite covering of $K$ such that by means of Lemma \ref{Lemma}, for fix index $i$, $V_i$ is homeomorphic to the open unit ball $B$ of ${\R}^n$. Let $(V_i ,\psi _i )$ be the corresponding chart, we complete the proof by means of partition of unity. More intuitively, let $\{\alpha _i\}$ be a partition of unity subordinate to the covering $\{V_i \}_{i=1}^m $. For approximating $\phi _j$ by $C^ \infty$ functions with compact support in $H_1^p (M)$, it remains to approximate each $\alpha _i \phi _j$ for $1 \leq i\leq m.$ For fixed $i,$ $\psi _i$ is a homeomorphic map between $V_i$ and the unit ball $B$. Consider the functions $(\alpha _i \phi _j)\circ \psi_i ^{-1}$ which have their support in $B$. Let us denote $u:=(\alpha _i \phi _j)\circ \psi_i ^{-1}$ to be consistent with notations of Lemmas \ref{lem1} and \ref{lem2}. Consider the convolution $J_\epsilon *u$ with $\lim\limits _{\epsilon \rightarrow o^+} J_\epsilon *u =u.$ Let $h'_\epsilon :=J_\epsilon *u \in C^\infty ({\R}^n),$ then $h'_\epsilon$ has a compact support, that is $h'_\epsilon \in D(B)$. We approximate $u$ by $h_k :=h'_\frac {\epsilon }{k}$. More precisely
$$\lim\limits _{k\rightarrow \infty } h_k = \lim\limits_{k\rightarrow \infty} h'_\frac{\epsilon}{k}=\lim\limits_{\epsilon \rightarrow o^+} h'_\epsilon=u= (\alpha _i \phi _j)\circ \psi_i ^{-1}.$$ Moreover,  $h_k$ is $C^\infty$.
Hence $h_k$ converges to $u$ in $H_1^p (B).$
Now $h_k \circ\psi _i$ converges to $\alpha _i \phi _j$ in $H_1^p (V_i)$ and $h_k \circ\psi _i \in D(V_i ).$
Thus we have approximated each $\alpha _i \phi _j$ by functions in $D(V_i )$. This completes the proof.
\hspace{\stretch{1}}$\Box$\\

 The similar proof can be repeated for backward geodesically complete spaces.

\begin{example}
Let $(M,F)$ be a Compact Finsler manifold. It is forward geodesically complete, hence by Theorem \ref{Th1}, $D(M)$ is dense in $H_1 ^p (M)$.
\end{example}

\bc\label{Cor;one}
 Let $(M,F)$ be a compact, connected, $C^\infty,$ reversible Finsler manifold and $f:M\longrightarrow{\R}$ a real function for which $\int _M fdv_F =0$, then the weak solution $u$ of the Dirichlet equation $\Delta u=f$ can be approximated by $C^\infty$ functions with compact support on $M$.
\ec

\begin{proof}\label{Ex;Dirichlet}
Theorem \ref{Th1} can be used to approximate weak solutions of Dirichlet problem on Finsler manifolds. Indeed in complete Finsler manifolds with certain conditions the Dirichlet problem $\Delta u=f$ has a unique solution which lies in Sobolev space $H_1^2 (M)$, for similar proof one can refer to \cite{Au2} and \cite{Ya}. Hence by Theorem \ref{Th1} we can approximate these weak solutions by $C^\infty$ functions with compact support on $M$.
\end{proof}

\subsection{Case of manifolds with $C^r$ boundary}
In proof of Theorem \ref{Th2} we use the technic applied in proof of the following theorem on half-spaces on ${\R}^n.$\\
\bt \cite{Au2}
$C^\infty (\overline E)$ is dense in $H^p_k (E)$, where E is a half-space $E=\{x \in {\R}^n : x_1 <0\}$ and $C^\infty (\overline E)$ is the set of functions that are restriction to $\overline E$ of $C^\infty$ functions on ${\R}^n.$
\et
\noindent {\bf Proof of Theorem \ref{Th2}. } Let $\phi$ be a real $C^\infty$ function on the Sobolev space $ H_k^p (W)$, that is, $\phi \in C^ \infty (W) \cap H_k^p (W). $
Here we approximate $\phi$ by the functions in $C^r (\overline W)$. Since $\overline W$ is compact, we can consider $(V_i,\psi _i ), i=1,\cdots ,N$ as a finite $C^r$ atlas on $\overline W.$ Each $V_i$ depending on $V_i \subset W$ or $V_i$ has intersection with $\partial W,$ is homeomorphic either to the unit ball of ${\R}^n$ or a half-ball $D=B\cap \overline E,$ respectively, where $E=\{x\in {\R}^n : x_1 <0\}.$ Let $\{\alpha _i\}$ be a $C^\infty$ partition of unity subordinate to the finite covering $\{V_i\}_{i=1}^m$ of $\overline W$. By properties of partition of unity,  it remains to show that each $\alpha _i \phi$, supported in $V_i,$ can be approximated by functions in $C^r (V_i).$ Each $V_i$ is homeomorphic either to the unit ball $B$ or a half-ball $D.$ First let $V_i$ be homeomorphic to $B$, then by the relation $V_i \subset W$ we have $\alpha _i \phi \in C^ {\infty} (V_i)$, therefore $\alpha _i \phi \in C^ r (V_i)$.\\
Now let $V_i$ be homeomorphic to $D=B\cap \overline E$ and $\psi _i$  a homeomorphism between $V_i$ and $D.$ Consider the sequence of functions $h_m$ as restricted to $D$ of $((\alpha _i \phi )\circ\psi_i^{-1})(x_1 -\frac{1}{m},x_2 ,x_3, \cdots ,x_s )$. Let $\phi \in H_k^p (W) \cap C^\infty (W),$ by appropriate choice of $\{V_i\}$ and $\{\alpha _i\}, $ the restriction of $((\alpha _i \phi )\circ\psi_i^{-1})(x_1 -\frac{1}{m},x_2 ,x_3, \cdots ,x_s )$ to $D$ and its derivative up to order $r$ converge to $((\alpha _i \phi )\circ\psi_i^{-1})$ in $H_k^p (D)$, where $D$ has the Euclidean metric. Therefore by Proposition \ref{TTT} we obtain
\begin{align*}
(\int_ {SM} |\nabla^t&(h_m\circ\psi_i-\alpha_i\phi)|^p dV_{SM})=\\
&\sum\limits_{i=1} ^{\infty}\int _{V_i} (\int _{S^{n-1}} |\nabla ^t (h_m\circ\psi_i-\alpha_i\phi)|^p _{(x,\frac{y}{F(x,y)})}\frac{det(g_{ij}(x,y))}{F(x,y)^n}d\sigma)dx\rightarrow 0,
\end{align*}
where $0\leq t \leq k$ is an integer.
Hence $h_m \circ\psi _i \longrightarrow \alpha _i \phi $ in $H_k^p (V_i),~k\leq r$, $h_m \circ\psi _i \in C^r (\overline W)$ and proof is complete.~ 
\hspace{\stretch{1}}$\Box$\\

In the following example we show that the assumption $k\leq r$ in Theorem \ref{Th2} is sharp and can not be omitted.
\begin{example}
Let $\overline W=[-1,1]\times[0,1]$ be a manifold with boundary of class $C^0$. Denote the points of $W$ by $x=(x^1 ,x^2)$ and the points of $T_x (W)$ by $y=(y^1 , y^2)$. Let $F(x,y)= \sqrt{g_x (y,y)}$ be a Finsler structure defined by
$g=g_{ij} ~dx^i \otimes dx^j=(dx^1)^2 + (dx^2)^2$ on $W$. Define, $u \in H_1^p (W)$ by $u:W\rightarrow {\R}$ and
\begin{equation}
u(x^1 ,x^2)=\left\{ \begin{array}{ll}1& \quad x^1> 0 \\ 0& \quad x^1\leq 0\end{array}.\right.
\end{equation}
 We claim that for sufficiently small $\varepsilon$, there is no $\phi \in C^1 (\overline W)$ such that $\parallel u-\phi\parallel_{H^p_1 (W)}<\varepsilon$. Assume for a while that our assumption is not true and the function $\phi$ exists. Let $S=\{(x^1 , x^2): -1 \leq x^1 \leq 0 ~ , 0 \leq x^2 \leq 1\}$ and $K= \{(x^1 , x^2) : 0<x^1 \leq 1 ~ , 0<x^2 \leq 1\}$, then $\overline W = S\cup K.$ On $S$, $u(x^1,x^2)=0$, hence $\parallel 0-\phi \parallel _{H_1^p (S)} < \varepsilon $ or $\parallel \phi \parallel _1 + \parallel \nabla\phi \parallel _1 < \varepsilon ,$ therefore $\parallel \phi \parallel _1 < \varepsilon$. On $K$, $u(x^1 ,x^2)=1$ thus $\parallel 1-\phi \parallel_1 < \varepsilon$ or $\parallel \phi \parallel_1 > 1-\varepsilon$. Put $\psi (x^1)= \int _0 ^1 \phi (x^1 , x^2)dx^2$, then there exist the real numbers $a$ and $b$ with $-1 \leq a \leq 0$ and $0<b \leq 1$ such that $\psi (a)<\varepsilon$ and $\psi (b)>1-\varepsilon$. Thus
\begin{align*}
1-2\varepsilon < \psi (b)- \psi (a) &= \int _a ^b \psi' (x^1)dx^1 \leq \int _{\overline W} |D_{x^1} \phi (x^1 , x^2)| dx^1 dx^2 \\
&\leq \big(\int _{\overline W} 1^{p'} dx^1 dx^2 \big)^ \frac {1}{p'} \big(\int _{\overline W} |D_{x^1} \phi (x^1 , x^2)|^p dx^1 dx^2 \big)^ \frac {1}{p}\\ &=2^{\frac{1}{p'}} \parallel D_{x^1} \phi (x^1 , x^2) \parallel _{L^p (\overline W)} <2^{\frac{1}{p'}} \varepsilon,
\end{align*}
where $\frac{1}{p} + \frac{1}{p'} =1$.
Hence $1<(2+2^{\frac {1}{p'}})\varepsilon$ which is not possible for small $\varepsilon$. This contradict our provisional assumption and prove the statement.
\end{example}
For some other Sobolev inequalities in Finsler geometry one can refer to \cite{Oh}.

{\small Faculty of Mathematics and computer science, Amirkabir University of Technology (Tehran Polytechnic),}\\
{\small 424 Hafez Ave. 15914, Tehran, Iran.}\\ {\small  bidabad@aut.ac.ir; alirezashahi@aut.ac.ir}
\end{document}